\documentclass[a4paper]{article}

\usepackage[english]{babel}
\usepackage[utf8x]{inputenc}
\usepackage[colorinlistoftodos]{todonotes}
\usepackage{amsmath,amsfonts,amssymb,hyperref,color,graphicx,siunitx,fancyhdr,lastpage,enumerate, commath, amsthm, subcaption}
\newtheorem{theorem}{Theorem}
\newtheorem{lemma}[theorem]{Lemma}
\newtheorem{proposition}[theorem]{Proposition}
\newtheorem{assumption}{Assumption}
\newtheorem{corollary}[theorem]{Corollary}

\theoremstyle{remark}
\newtheorem{remark}[theorem]{Remark}

\numberwithin{theorem}{section}

\usepackage{url}
\usepackage[numbers]{natbib}

\newcommand{\R}{\mathbb R}
\newcommand{\Rd}{\mathbb R^d}
\newcommand{\E}{\mathbb E}
\renewcommand{\P}{\mathbb P}
\renewcommand{\Pr}{\mathcal P(M)}

\newcommand{\Lk}{L^\kappa}

\newcommand{\taud}{\tau_\partial}
\newcommand{\lkap}{\lambda_0^\kappa}
\newcommand{\Res}{ \mathcal R}
\newcommand{\etR}{e^{t\Res}}
\newcommand{\Kt}{\mathcal K^t}
\newcommand{\ubar}[1]{\text{\b{$#1$}}}

\newcommand{\footremember}[2]{%
    \footnote{#2}
    \newcounter{#1}
    \setcounter{#1}{\value{footnote}}%
}
\newcommand{\footrecall}[1]{%
    \footnotemark[\value{#1}]%
}

\title{An approximation scheme for quasi-stationary distributions of killed diffusions}
\author{%
Andi Q. Wang\footremember{Oxf}{University of Oxford}\footnote{Corresponding author. Address: Department of Statistics, 24--29 St Giles', Oxford, OX1 3LB, UK. Email: \href{mailto:a.wang@stats.ox.ac.uk}{a.wang@stats.ox.ac.uk}.}
\and Gareth O. Roberts\footnote{University of Warwick}
\and David Steinsaltz\footrecall{Oxf}}

\date{July 1, 2019}

\begin{document}
\maketitle

\begin{abstract} 
In this paper we study the asymptotic behavior of the normalized weighted empirical occupation measures of a diffusion process on a compact manifold which is killed at a smooth rate and then regenerated at a random location, distributed according to the weighted empirical occupation measure. We show that the weighted occupation measures almost surely comprise an asymptotic pseudo-trajectory for a certain deterministic measure-valued semiflow, after suitably rescaling the time, and that with probability one they converge to the quasi-stationary distribution of the killed diffusion. These results provide theoretical justification for a scalable quasi-stationary Monte Carlo method for sampling from Bayesian posterior distributions.

\sloppy \textit{Keywords:} asymptotic pseudo-trajectory, killed diffusion, quasi-stationary distribution, quasi-stationary Monte Carlo method, stochastic approximation.
\end{abstract}

\section{Introduction}
In this paper we consider a stochastic approximation algorithm to approximate the quasi-stationary distribution of a killed diffusion on a compact manifold.
This is particularly motivated by recent developments in theoretical and computational statistics: a promising new paradigm has emerged for the performance of exact Bayesian inference on large datasets, known as \textit{quasi-stationary Monte Carlo} (QSMC) methods, \cite{Pollock2016, Wang2017}. Recall that for a Markov process $Y=(Y_t)_{t\ge 0}$ killed at random time $\taud$, a distribution $\pi$ is called \textit{quasi-stationary} if for each $t\ge 0,$
\begin{equation}
\P_\pi (Y_t \in \cdot \,|\tau_\partial >t) = \pi(\cdot)
\label{eq:QSD}
\end{equation}
where $\P_\pi = \int \pi(\dif y) \,\P_y$, $\P_y$ denoting the law of $Y$ with initial position $Y_0=y$. The idea behind quasi-stationary Monte Carlo is to draw samples from such a quasi-stationary distribution. 

There are (at least) two distinct approaches to realising this method in practice. The first is to use particle approximation methods. This was the approach of the original QSMC method, the Scalable Langevin Exact (ScaLE) algorithm, which was introduced in \cite{Pollock2016}. The quasi-stationary framework enables the principled use of \textit{subsampling} techniques, and the resulting ScaLE algorithm is provably efficient for performing exact Bayesian inference when the underlying dataset is large. This approach is supported by the full weight of the sequential Monte Carlo literature, and thus has theoretical assurances of convergence.

While the ScaLE algorithm possesses some desirable theoretical and computational properties, it also suffers several drawbacks. One particularly prominent drawback is the high complexity of implementing the algorithm. A cursory glance through \cite{Pollock2016} can convince the reader that the algorithm is (necessarily) involved. This motivates the search for alternative QSMC methods which inherit the desirable properties of exactness and scalability, while avoiding the algorithmically complex particle approximation framework. 

One alternative approach to QSMC takes advantage of the fact that quasi-stationary distributions can be written as solutions to fixed point equations in measure spaces (represented here by equation \eqref{eq:fixed_pt} in Proposition \ref{prop:QSD_exist}). 
This yields the QSMC method dubbed Regenerating ScaLE (ReScaLE), which seeks to approximate the quasi-stationary distribution $\pi$ by means of stochastic approximation. The idea behind ReScaLE is that we simulate a \textit{single} killed diffusion path, whose quasi-stationary distribution coincides with the Bayesian posterior of interest. When this single particle is killed, the particle is `regenerated' or `reborn' in a new random location drawn independently from its normalized weighted empirical measure. It then continues to evolve from this new starting position until it is killed again, at which point it is again regenerated, and so on. This algorithm is inspired by similar algorithms for the approximation of quasi-stationary distributions in discrete-time settings \cite{Aldous1988, Benaim2015, Blanchet2016, Benaim2016}. Indeed, even outside of QSMC such algorithms have been applied to the simulation of quasi-stationary distributions in ecology \cite{Schreiber2018}.

The resulting ReScaLE algorithm inherits the key properties that motivate the ScaLE algorithm, and algorithmically it is significantly more transparent and straightforward, but this is worthwhile only if
the method provably leads to correct results. 
Theoretical analysis of stochastic approximation approaches to numerically solving fixed point equations are well-studied in finite-dimensional contexts, {\em cf.} \cite{Kushner2003}, but there is currently no theory appropriate for the measure-valued, continuous-time context of ReScaLE. While the implementation of the original ScaLE algorithm is underpinned by the vast body of work on sequential Monte Carlo, there has been no theoretical assurance that the ReScaLE algorithm even converges. The purpose of the present work
is to demonstrate this fundamental convergence property, leaving practical and computational properties of the algorithm to be considered in detail in forthcoming work.


In this work we assume that we are working on a compact, boundaryless, connected, $d$-dimensional smooth Riemannian manifold $M$, as in \cite{Benaim2002}. Relaxing the assumption of compactness is discussed briefly in Section \ref{sec:extens}. We have a particle $X=(X_t)_{t\ge 0}$ evolving on $M$ in continuous time, according to the solution to the stochastic differential equation (SDE)
\begin{equation}
\dif Y_t = \nabla A(Y_t)\dif t+\dif W_t
\label{eq:SDE}
\end{equation}
between regeneration events,
where $A:M\to \R$ is a smooth function and $W$ is a standard Brownian motion on $M$. Regeneration events occur at a state-dependent rate $\kappa(X_{t-})$, where $\kappa:M\to [0,\infty)$ is a given smooth, positive function, which we will refer to as the \textit{killing rate}. At a regeneration event, say at time $T$, the particle is instantaneously `killed' and `reborn': its location is drawn (independently) from its normalized weighted empirical occupation measure $\mu_T$, where for all $t\ge 0$, $\mu_t$ is given by
\begin{equation*}
\mu_t = \frac{r \mu_0}{r+\int_0^t \eta_s \dif s}  +\frac{\int_0^t \eta_s\delta_{X_{s-}}\dif s}{r+\int_0^t \eta_s \dif s} 
\label{eq:mu_t}
\end{equation*}
where $r>0$, $\eta_\cdot: \R_+ \to \R_+$ and a probability measure $\mu_0$ on $M$ are fixed. The resulting process $X$ is clearly non-Markovian. The addition of the $\mu_0$ term has the benefit of regularising the $\mu_t$ around $t=0$, as well as providing practical flexibility for the resulting Monte Carlo algorithm. The weight function $\eta_\cdot$ similarly provides additional practical flexibility; a straightforward choice would be constant $\eta_s \equiv 1$, but see Remark \ref{rmk:eta_choice} for a nonconstant alternative.

\subsection{Main results}
The goal of this paper is to characterize the asymptotic behavior of the measure-valued
process $(\mu_t)_{t\ge 0}$, and to show that it converges to the quasi-stationary distribution
for the original killed process. We proceed in two steps, following the ``ODE method'' (\textit{cf.} \cite{Benaim1999, Kushner2003}), which has been used to prove convergence of similar reinforced processes, for instance in \cite{Benaim2015,Benaim2016, Benaim2002, Kurtzmann2010}. The ODE method proceeds by two key steps. First,  showing that a certain deterministic semiflow $\Phi$ converges to the appropriate limit, in our case the quasi-stationary distribution $\pi$. Second,
showing that, following a suitable deterministic time change $\zeta_t := \mu_{h(t)}$, the stochastic evolution of the measures $(\zeta_t)_{t\ge 0}$ shadows $\Phi$
in an appropriate sense, to be defined below. From these two properties 
almost-sure convergence of $\mu_t$ to the quasi-stationary distribution can be deduced.

The present paper extends previous related work in considering a 
continuous-time diffusive process on a compact manifold (rather than, say, a discrete-time Markov chain), which experiences `soft killing' according to a smooth state-dependent killing rate (rather than instantaneous `hard killing' at a boundary). In this setting of soft killing, a generic killing time $\taud$ has the form
\begin{equation} 
    \taud = \inf\left\{t\ge 0: \int_0^t \kappa(X_{s})\dif s > \xi\right\},
    \label{eq:taud}
\end{equation}
where $\xi\sim\text{Exp}(1)$ is independent of the process $X$, rather than the first hitting time of some forbidden set of states as in the hard killing case.

\begin{theorem}
Under Assumptions \ref{assump:M}, \ref{assump:A_smooth}, \ref{assum:kappa} and \ref{assump:eta} detailed in Section \ref{subsec:assumps}, we have with probability 1, that for each $T>0$,
\begin{equation}
\lim_{t\rightarrow \infty} \sup_{ s\in[0,T]} d_w\big (\zeta_{t+s}, \Phi_s(\zeta_t)\big) = 0
\label{eq:APT}
\end{equation}
where $d_w$ is a metric that metrises weak-* convergence of probability measures given in \eqref{eq:metric} and $\Phi$ is the semiflow defined in Section \ref{sec:flow}.
\label{thm:intro_apt}
\end{theorem}

\begin{remark}
In \cite{Benaim1999}, a map $t\mapsto \zeta_t$ that satisfies \eqref{eq:APT} for each $T>0$ is called an \textit{asymptotic pseudo-trajectory} for the measure-valued semiflow $\Phi$.
\end{remark}

In Section \ref{sec:flow} we define the semiflow $\Phi$, and proceed to show that it has a global attractor $\pi$, which is the unique quasi-stationary distribution of the diffusion \eqref{eq:SDE} killed at rate $\kappa$.

In particular, Theorem \ref{thm:intro_apt} leads to the following corollary:
\begin{corollary}
Under the conditions of Theorem \ref{thm:intro_apt}, we have that almost surely, $\lim_{t\to\infty} \mu_t = \pi$ in the sense of weak-* convergence.
\label{cor:conv}
\end{corollary}
\noindent That is, we have $\lim_{t\to\infty}\mu_t(f) = \pi(f)$ for any continuous $f:M\to \R$.

\begin{figure}

\begin{subfigure}{.5\textwidth}
  \centering
  \includegraphics[width=.65\linewidth]{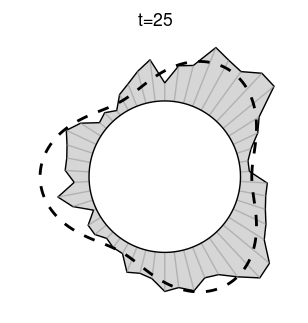}
\end{subfigure}
\begin{subfigure}{.5\textwidth}
  \centering
  \includegraphics[width=.65\linewidth]{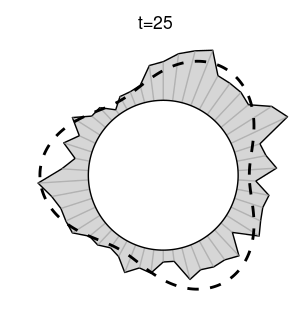}
\end{subfigure}
\begin{subfigure}{.5\textwidth}
  \centering
  \includegraphics[width=.65\linewidth]{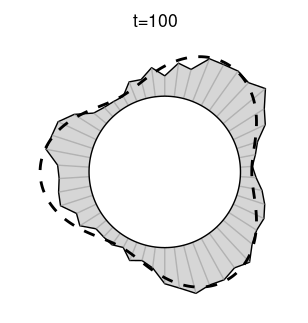}
\end{subfigure}
\begin{subfigure}{.5\textwidth}
  \centering
  \includegraphics[width=.65\linewidth]{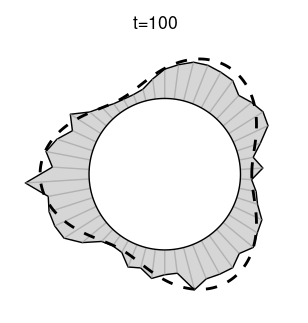}
\end{subfigure}
\begin{subfigure}{.5\textwidth}
  \centering
  \includegraphics[width=.65\linewidth]{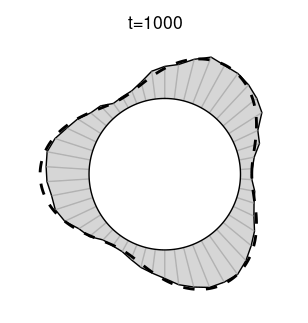}
\end{subfigure}
\begin{subfigure}{.5\textwidth}
  \centering
  \includegraphics[width=.65\linewidth]{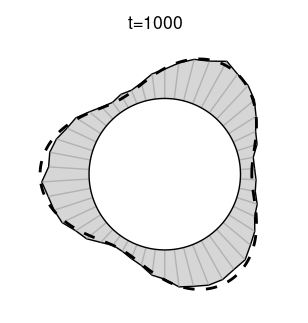}
\end{subfigure}
\begin{subfigure}{.5\textwidth}
  \centering
  \includegraphics[width=.65\linewidth]{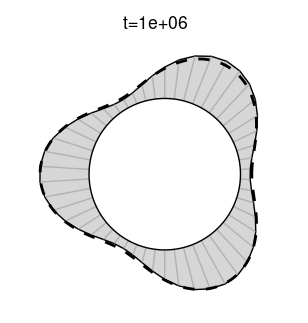}
\end{subfigure}
\begin{subfigure}{.5\textwidth}
  \centering
  \includegraphics[width=.65\linewidth]{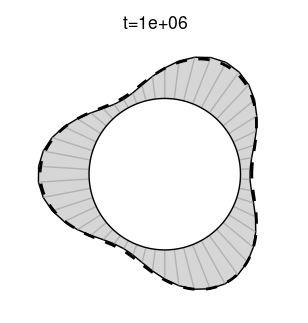}
\end{subfigure}

\caption{Example on the circle. The circle is parameterised by $\theta \in [0,2\pi)$, with due east corresponding to $\theta=0$. The underlying diffusion is a Brownian motion, and we have chosen $r=1000$, $\mu_0$ to be the uniform distribution on $M$ and $\eta_s \equiv 1$. The quasi-stationary distribution $\pi(\theta)=(0.3+ \sin^2(1.5\theta))/(1.6\pi)$ is the dashed line. We have plotted $\int_0^t \delta_{X_s-}\dif s/t$ discretised into 50 evenly-spaced bins, for $t=25, 100, 1000, 10^6$. The two columns are two independent runs.}
\label{fig:circ_dens}
\end{figure}

While our results above hold for any appropriate given killing rate which satisfies our assumptions, in the ReScaLE algorithm the killing rate $\kappa$ is chosen so that the quasi-stationary distribution $\pi$ \textit{equals} the Bayesian posterior distribution of interest; see the expression \eqref{eq:kappa}, \textit{cf.} \cite{Wang2017}. Corollary \ref{cor:conv} tells us that in this setting we can draw approximate samples from $\pi$ by running the regenerating process $X$ and outputting $\mu_t$ for a large $t$ as a proxy for $\pi$.

Figure \ref{fig:circ_dens} shows the output of two independent simulations on the unit circle $M=\mathbb R/(2\pi\mathbb Z)$ parameterised by $\theta\in[0,2\pi)$, which is amenable to straightforward simulation and visualisation. The underlying diffusion is a Brownian motion ($A\equiv 0$ in \eqref{eq:SDE}), and the quasi-stationary distribution is trimodal with density function $\pi(\theta)=(0.3+ \sin^2(1.5\theta))/(1.6\pi)$ with respect to the Riemannian measure. The corresponding killing rate, calculated using \eqref{eq:kappa} is $\kappa(\theta)=2.25(2\cos^2(1.5\theta)-1)/(0.3+\sin^2(1.5\theta))+K$, with $K=1.75$.

The simulations were done using a simple Euler discretisation scheme where time was discretised into intervals of length 0.05. The plots are oriented so that $\theta =0$ is due east, with $\theta=\pi/2$ being due north, and so on. We have chosen $r=1000$, $\mu_0$ to be the uniform distribution over $M$ and $\eta_s \equiv 1$. The plots shows $\int_0^t \delta_{X_s-}\dif s/t$ for $t=25, 100, 1000, 10^6$, split into 50 evenly-spaced bins. The quasi-stationary distribution $\pi$ is the dashed line.

We see that there is a significant amount of variability for the small time values, as the initial rebirths are drawn mostly from the uniform distribution $\mu_0$. These discrepancies are largely smoothed out by $t=1000$ and certainly by $t=10^6$, $\int_0^t \delta_{X_s-}\dif s/t$ closely resembles the quasi-stationary distribution $\pi$.

This present work closely follows the approaches of \cite{Benaim2002} and the recent \cite{Benaim2016} in order to prove our limiting result Theorem \ref{thm:intro_apt}. \cite{Benaim2002} shows an analogous characterization for the normalized occupation measures of a self-interacting diffusion on a compact space, where the empirical occupation measure influences the present behavior of the diffusion through its drift term. In our work, the influence of the occupation measures is felt through jumps, which occur at a state-dependent rate. In \cite{Benaim2016}, the authors prove a discrete-time analogue of our above result; their underlying Markov process is a discrete-time Markov chain evolving on a compact space, rather than a diffusion. In \cite[Section 8.3]{Benaim2016}, the authors do suggest a continuous-time extension of their work: a diffusive process that is killed instantaneously when hitting the boundary of an open set. This is not the same as the present work;  we are assuming a boundaryless manifold, and instead of hard killing at a boundary, killing occurs at a smooth state-dependent rate $\kappa$ as in \eqref{eq:taud}. The key difference in the proposed setting of \cite[Section 8.3]{Benaim2016} is that in their case, the state space is no longer compact, and hence additional arguments ensuring almost sure tightness of the empirical occupation measures are needed.


We follow generally the path mapped out by \cite{Benaim2002, Benaim2016, Kurtzmann2010}, often referred to as the ``ODE method'', \textit{cf.} \cite{Benaim1999, Kushner2003}. To understand the particular continuous-time dynamics we employ techniques similar to those of \cite{Wang2017, Kolb2012}. For example, to handle the killing mechanism we make use of the transition subdensity of the killed diffusion and the corresponding resolvent operator; see Lemmas \ref{lemma:subtrans}, \ref{lem:Res_basic} and \ref{lem:inv_pi}. The continuous-time setting also provides a natural interpretation of the weights $\eta_s$ in terms of the distribution of the rebirth `times'; see Remark \ref{rmk:eta_choice}. Another contribution of this work is the analysis of the deterministic semiflow in Section \ref{sec:flow}, where we transform to an auxiliary Markov process and apply an appropriate drift condition.  

The earliest version of a similar analysis that we are aware of is the convergence
proof of \cite{Aldous1988} for the finite-state-space discrete-time setting.
\cite{Benaim2015} and \cite{Blanchet2016} expanded this work to derive rates of convergence and central limit theorems. These analyses rely upon special
techniques, such as those in \cite{Kushner2003}, applicable only to finite-dimensional
probability distributions.

The structure of the paper is as follows. In Section \ref{sec:prelim} we lay out the mathematical setting and describe the assumptions that we are making. We also define the fixed rebirth process and prove some of its key properties which are crucial for later defining the semiflow. In Section \ref{sec:flow} we define and analyse the deterministic semiflow $\Phi$ and prove that $\pi$ is a global attractor. In Section \ref{sec:APTs} we prove that our normalized weighted occupation measures almost surely comprise an asymptotic pseudo-trajectory for $\Phi$, concluding the convergence proof. Finally, in Section \ref{sec:extens}, we briefly discuss the possibilities for extending these results
to noncompact manifolds.

\section{Preliminaries}
\label{sec:prelim}
\subsection{General background}
\label{subsec:assumps}
We first describe our assumptions and some notation. These assumptions are \textit{assumed to hold throughout the paper}.
\begin{assumption}
    $M$ is a $d$-dimensional boundaryless $C^\infty$ compact connected Riemannian manifold.
    \label{assump:M}
\end{assumption}
We will denote the corresponding Riemannian measure by $\dif x$ or $\dif y$, and the Riemannian inner product of points $x$ and $y$ by $x\cdot y$. Let $C(M)$ denote the Banach space of (bounded) real-valued continuous functions on $M$ equipped with the sup norm, $\|\cdot\|_\infty$. Let $\Pr$ denote the space of Borel probability measures on $M$, equipped with the topology of weak-* convergence, namely convergence along all bounded continuous test functions. (This is conventionally called simply weak convergence in probability theory, but we follow the terminology of \cite{Benaim2016} to avoid confusion when working with the Banach space $C(M)$ and its dual.) Thus $\Pr$ is a compact metrisable space, by the Prokhorov theorem, {\it cf.} \cite[Theorem 11.5.4]{Dudley2002}. 

As $M$ is compact, $C(M)$ is separable, so we may choose a sequence of smooth functions $f_1, f_2, \dots$ which are dense in $\{f\in C(M): \|f\|_\infty\le 1\}$.
We define the metric
$d_w:\Pr \times \Pr \to \R_+$ on $\Pr$ by
\begin{equation}
d_w(\nu_1, \nu_2) = \sum_{i=1}^\infty \frac{1}{2^i}|\nu_1(f_i)-\nu_2(f_i)|.
\label{eq:metric}
\end{equation}
We also define the total-variation norm for a signed measure $\mu$ on $M$ by
\begin{equation*}
\|\mu\|_1 = \sup\{|\mu(f)|: f\in C(M), \|f\|_\infty \le 1 \}.
\end{equation*}

Let $\Omega$ be the Skorokhod space of c\`{a}dl\`{a}g paths $\omega: \R_+ \to M$, and let $\mathcal F$ be the cylinder $\sigma$-algebra. Let $X=(X_t)_{t\ge 0}$ be the coordinate process, $X_t(\omega)=\omega(t)$, and let $(\mathcal F_t)_{t\ge 0}$ be the natural filtration of $X$.

\begin{assumption}
    The function $A: M \to \R$ is $C^\infty$.
    \label{assump:A_smooth}
\end{assumption}
In particular, this implies that $A$ is bounded (since $M$ is compact), and that the SDE \eqref{eq:SDE} on $M$ has a unique strong solution for any initial position.

Recall that for a Markov process $Y$ killed at a random stopping time $\tau_\partial$, a distribution $\pi\in \Pr$ is \textit{quasi-stationary} if it satisfies \eqref{eq:QSD} for all $t\ge 0$.

Given a measurable killing \textit{rate} $\kappa:M\to [0,\infty)$, the corresponding killing time $\tau_\partial$ is defined as
\begin{equation*}
\tau_\partial := \inf\bigg\{ t\ge 0: \int_0^t \kappa(Y_s)\dif s\ge  \xi \bigg \},
\end{equation*}
where $\xi$ is an independent exponential random variable with rate 1.

\begin{assumption}
The killing rate $\kappa:M\to [0,\infty)$ is $C^\infty$ and is uniformly bounded away from zero: there exists some constant $\ubar{\kappa}>0$ such that 
\begin{equation}
0<\ubar{\kappa}\le \kappa(x)  \quad \forall x \in M.
\label{eq:lower_bd_kap}
\end{equation}
As a continuous function on a compact space $\kappa$ is necessarily bounded above, say $$\kappa(x)\le \bar\kappa<\infty\quad \forall x\in M.$$
\label{assum:kappa}
\end{assumption}
Given a killing rate $\kappa$ which is not strictly bounded away from 0, we can always add a positive constant everywhere; this will not affect the quasi-stationary behavior of the process and will ensure that \eqref{eq:lower_bd_kap} holds. The upper bound on the killing rate will certainly guarantee that the resulting process will be almost surely nonexplosive.

We will later see that given $\kappa$, the diffusion \eqref{eq:SDE} killed at rate $\kappa$ has a unique quasi-stationary distribution (Proposition \ref{prop:QSD_exist}).

The question of existence and uniqueness of quasi-stationary distributions of killed diffusions has been studied in depth, for instance in \cite{Kolb2012}, \cite[Chapter 6]{Collet2013}. This context of quasi-stationary Monte Carlo methods -- where we are starting from a density $\pi$ and wish to construct a killed diffusion whose quasi-stationary distribution coincides with $\pi$ -- is the topic of \cite{Wang2017}, and the interested reader is encouraged to look therein for more details.

From \cite{Wang2017}, if we were to start with a smooth positive density $\pi$, the appropriate choice of $\kappa$ is given for each $y\in M$ by
\begin{equation}
 \kappa (y) := \frac{1}{2}\bigg (\frac{\Delta \pi}{\pi} - \frac{2\nabla A \cdot \nabla \pi}{\pi} - 2\Delta A \bigg)(y)+K.
\label{eq:kappa}
\end{equation}
Here $\nabla$ and $\Delta$ are the gradient and Laplacian operators on $M$, and $K$ is a constant chosen so that $\kappa$ satisfies \eqref{eq:lower_bd_kap}. This choice of $\kappa$ is discussed at length in \cite{Wang2017}. In this work we will not necessarily assume we are starting from $\pi$, and take $\kappa$ to be a general killing rate satisfying Assumption \ref{assum:kappa}. 

Fix $\mu_0 \in \Pr$ and $r>0$. We fix a weight function $\eta:\R_+ \to \R_+$ satisfying the following assumptions: Define functions $g:\R_+\to \R_+$ and $\alpha:\R_+ \to \R_+$ by
\begin{equation*}
g(t) := \int_0^t \eta_s \dif s,\quad \alpha_t := \frac{\eta_t}{r+g(t)}, \quad t \ge0.
\end{equation*}

\begin{assumption}
    $\eta$ is continuously differentiable with $\eta_t>0$ for all $t>0$, and $g(t)\to\infty$ as $t\to \infty$. $\alpha$ is differentiable and satisfies $\alpha_t \to 0$ as $t\to\infty$, $\int_0^\infty \alpha_s^2 \dif s<\infty$ and  
    \begin{equation}
        \sum_{n=1}^\infty \int_{h(n)}^\infty \alpha_s^2 \dif s <\infty,
        \label{eq:gamma_int_sum}
    \end{equation}
    where $h$ is defined below.
    \label{assump:eta}
\end{assumption}
Since $g$ is strictly increasing (as $\eta_t >0$), continuously differentiable and increases to $\infty$, it is a diffeomorphism of $\R_+ \to \R_+$. Thus it has a well-defined continuously differentiable inverse $g^{-1}$. The function $h:\R_+ \to \R_+$ is then given by $$h(t):=g^{-1}(r\mathrm e^t-r), \quad t \ge 0.$$
This function $h$ will be the time change that we shall employ in Section \ref{sec:APTs}.

\begin{remark}
It follows from Assumption \ref{assump:eta} that $\int_0^t \alpha_s\dif s = \log(1+g(t)/r)\to\infty$ as $t\to \infty$. Thus these conditions on $\alpha$ are analogous to the typical discrete-time assumptions on the step sizes in traditional stochastic approximation; \textit{cf.} \cite[Chapter 5]{Kushner2003}. Since $t\mapsto \int_{h(t)}^\infty \alpha_s^2 \dif s$ is monotone decreasing, a sufficient condition for \eqref{eq:gamma_int_sum} to hold is that
\begin{equation*}
\int_0^\infty \dif t \int_{h(t)}^\infty \dif s \,\alpha^2_s =\int_0^\infty \alpha^2_s \log(1+g(s)/r)\dif s <\infty.
\end{equation*}
\end{remark}

\noindent Define the normalized empirical occupation measures $(\mu_t(\omega))_{t\ge 0}$ by 
\begin{equation*}
\mu_t(\omega) = \frac{r\mu_0+\int_0^t \eta_s \delta_{\omega(s-)}\dif s}{r+\int_0^t \eta_s\dif s}
\end{equation*}
where $\int_0^t \eta_s\delta_{\omega(s-)}\dif s(A) = \int_0^t \eta_s 1_A(\omega(s-))\dif s$ for each measurable $A\subset M$. In general we will omit the dependence on $\omega$.

\begin{remark}
For simplicity one may take $\eta_t \equiv 1$, then $g(t)=t$ and $\alpha_t = 1/(r+t)$ for all $t\ge 0$. Sampling from $Z\sim \int_0^t \eta_s \delta_{X_{s-}}\dif s/\int_0^t \eta_s \dif s$ is then equivalent to simulating $V\sim \text{Unif}([0,1])$ and setting $Z= X_{(Vt)-}$. More generally, for $k\ge 0$ one can take $\eta_t = t^k$ for each $t\ge 0$. It is not difficult to check that for this choice of $\eta$ Assumption \ref{assump:eta} is satisfied, and simulating $Z\sim \int_0^t \eta_s \delta_{X_{s-}}\dif s/\int_0^t \eta_s \dif s$ is equivalent to the case of constant $\eta$ except with $V\sim $ Beta$(k+1,1)$. Heuristically, choosing a larger value of $k$ prioritizes the more recent times. 
The choice of the parameter $k$ to accelerate convergence is itself an interesting question, which will be explored in future work. Preliminary simulations involving Brownian motions and unimodal targets seem to suggest $k=10$ might be a reasonable choice.
\label{rmk:eta_choice}
\end{remark}

For any $\mu\in\Pr$, define the operator $-L_\mu$ on smooth functions by
\begin{equation}
-L_\mu f(x) = \frac{1}{2}\Delta f(x) + \nabla A \cdot \nabla f(x)+\kappa(x) \int_M \big (f(y)-f(x)\big)\,\mu(\dif y).
\label{eq:gen_FRmu}
\end{equation}
Here we choose to use the negative operator in order to comport with the convention 
adopted in \cite{Kolb2012,Wang2017}, where it was chosen to make 
the corresponding self-adjoint operators positive. Formally, we are defining a generator on a core of smooth functions, which can be extended to a formal domain in the
standard way. However since we are working on a compact boundaryless manifold these technical
details raise no relevant complications, and hence are omitted.

We can
define probability measures $(\P_x:x\in M)$ with the following properties:
\begin{itemize}
\item $\P_x(X_0=x)=1$.
\item For all smooth $f\in C(M)$
\begin{equation*}
N_t^f := f(X_t) - f(x) -\int_0^t (-L_{\mu_s}) f(X_{s-})\dif s
\end{equation*}
is a $\P_x$-martingale with respect to $(\mathcal F_t)$.
\end{itemize}
We can construct $\P_x$ by explicitly defining the killing events, as follows. Let $\xi_1, \xi_2, \dots$ be an i.i.d. sequence of Exp(1) random variables. Set $T_0 = 0$, $X_0 =x$ and given $T_n$ and $X_{T_n}$ inductively define
\begin{equation}
\tau_{n+1} := \inf\bigg\{ t\ge 0: \int_0^t \kappa\big (Y_s^{(n+1)}\big)\dif s\ge \xi_{n+1} \bigg \}
\label{eq:killing_time}
\end{equation}
where for each $n=0,1,2,\dots$, $(Y_s^{(n+1)})_{s\ge 0}$ is an independent realisation of the solution to the SDE \eqref{eq:SDE} started from $Y_0^{(n+1)}=X_{T_n}$. For a careful treatment of defining diffusions on manifolds, the reader is referred to \cite{Stroock2000}.

Then update the path $X_{T_n +s} = Y^{(n+1)}_s$ for $s\in [0,\tau_{n+1})$. Set $T_{n+1}=T_n+\tau_{n+1}$ and independently draw $X_{T_{n+1}}\sim \mu_{T_{n+1}}$, which only depends on the path before time $T_{n+1}$.

\begin{remark}
    The practical simulation of this process, namely of the SDE dynamics \eqref{eq:SDE} and the killing times \eqref{eq:killing_time}, may seem difficult at first glance. For the killing times, in the present compact setting simulation is actually straightforward, since one can directly employ Poisson thinning (see \cite[Chater VI.2.4]{Devroye1986}) as $\kappa$ is bounded. Even in noncompact spaces simulation of such times may be performed without error through \textit{layered processes}; this is the case for the ScaLE algorithm \cite{Pollock2016}. To simulate the SDE, it turns out that through the techniques of \textit{exact simulation} \cite{Beskos2006}, in many settings the SDE \eqref{eq:SDE} can be simulated on fixed time horizons without error and without resorting to discretization. A thorough computational analysis of the resulting ReScaLE algorithm, including selection of the weights $\eta$ (\textit{c.f.} Remark \ref{rmk:eta_choice}), applications to large data sets and comparisons to existing methods will be the subject of a forthcoming piece of work.
\end{remark}

\subsection{Fixed Rebirth processes}
We now define the fixed rebirth processes and derive some useful properties. These will be crucial later for defining the deterministic semiflow. It will be convenient to work on $M$ with the measure $$\Gamma(\dif y)=\gamma(y)\dif y,$$ where $$\gamma(y) = \exp(2A(y)).$$

Let us write $\Lk$ for minus the generator of the diffusion $Y$ from \eqref{eq:SDE} killed at rate $\kappa$ and let $\mathcal D(L^\kappa)\subset C(M)$ be its domain. Such $f$ are twice continuously differentiable, and we have
\begin{equation*}
-\Lk f = \frac{1}{2}\Delta f  +\nabla A \cdot \nabla f - \kappa f.
\end{equation*}
Then we have the identity
\begin{equation}
\mathrm e^{-t\Lk}f(x) = \E_x[f(Y_t)1_{\{\tau_\partial>t\}}]=  \E_x \left [f(Y_t)\mathrm e^{-\int_0^t \kappa(Y_s)\dif s}\right ].
\label{eq:semigp_alt}
\end{equation}
Here (and throughout) $\tau_\partial$ denotes a general killing time, defined analogously to \eqref{eq:killing_time} and $Y$ evolves according to our SDE \eqref{eq:SDE} without any killing. The exponentiation of the operator $-L^\kappa$ is rigourously justified through the spectral theorem for self-adjoint operators, and the resulting (sub-Markovian) semigroup $(\mathrm e^{-tL^\kappa})_{t\ge 0}$ is strongly continuous on $C(M)$. Details may be found in \cite{Demuth2000}.

Similarly to \cite{Wang2017, Kolb2012}, we show the existence of a continous, positive transition subdensity for the SDE \eqref{eq:SDE} killed at rate $\kappa$.

\begin{lemma}
    The SDE \eqref{eq:SDE} killed at rate $\kappa$ has a $C^\infty$ positive transition subdensity $p^\kappa (t,x,y)$ with respect to $\Gamma$, that is,
    \begin{equation*}
        \mathrm e^{-t\Lk}f(x) = \int_M f(y)\, p^\kappa(t,x,y) \Gamma(\dif y).
    \end{equation*}
    \label{lemma:subtrans}
\end{lemma}
\noindent \textbf{Proof.}
First, we note that the \textit{unkilled} diffusion has a smooth positive transition density $p^0(t,x,y)$ with respect to $\Gamma$. The existence of this density is described briefly in Example 9 of \cite[Chapter 1.C]{Demuth2000}, and the details can be found in \cite[Chapter II]{Bismut1984}. In particular, the assumptions of \cite{Bismut1984} are that the manifold $M$ is a $C^\infty$ compact connected finite-dimensional Riemannian manifold (as in our Assumption \ref{assump:M}) and that the drift -- our $\nabla A$ -- is a $C^\infty$ vector field, which we are assuming (Assumption \ref{assump:A_smooth}).

In order to obtain the transition subdensity for the killed diffusion, we make use of \eqref{eq:semigp_alt}.
By conditioning on the end point we can write
\begin{equation*}
    \mathrm e^{-t\Lk}f(x) = \E_x\left[ f(Y_t) g(t,x,Y_t) \right]= \int_M f(y) g(t,x,y) p^0 (t,x,y)\Gamma(\dif y),
\end{equation*}
where for each $t>0$, $x,y\in M$, $$g(t,x,y) := \E_x\left[ \mathrm e^{-\int_0^t \kappa(Y_s)\dif s} \big | Y_t = y \right].$$
By \cite[Chapter II.d]{Bismut1984}, since $\kappa$ is bounded $C^\infty$, $g(t,x,y)$ is continuous in $x,y $ and in fact since $\kappa$ is smooth and nonnegative, $g(t,x,y)$ is jointly continuous over $t,x,y$ and smooth as a function of $x$ or $y$. Since $\kappa$ is nonnegative and bounded above, we have the bounds $ \mathrm e^{-t\bar\kappa}\le g(t,x,y)\le 1$. In particular $g(t,x,y)$ is positive, and hence setting for each $t>0$, $x,y\in M$,
\begin{equation}
    p^\kappa(t,x,y) = p^0(t,x,y) g(t,x,y),
    \label{eq:p^kappa}
\end{equation}
we obtain the positive  $C^\infty$ transition subdensity $p^\kappa(t,x,y)$ with respect to $\Gamma$ of diffusion \eqref{eq:SDE} killed at rate $\kappa$.  $\Box$


\begin{lemma}
The resolvent operator $\Res:C(M)\to C(M)$ given by
\begin{equation*}
\Res f(x) = \int_0^\infty \dif t \int_M \Gamma(\dif y) \,p^\kappa(t,x,y) f(y)
\end{equation*}
is a well-defined bounded, positive linear operator.
\label{lem:Res_basic}
\end{lemma}

\begin{remark}
We have defined $\Res$ as an operator on the Banach space $C(M)$. Its dual operator acts on the space of finite signed Borel measures on $M$ ({\em cf.} \cite[Section 7.4]{Dudley2002}). Following the standard probabilistic notation we will denote its dual action on a measure $\mu$ by simply $\mu \Res$. That is, $\mu \Res$ is the measure defined by
\begin{equation*}
\mu \Res (f) = \int_M \mu(\dif x ) \Res f(x).
\end{equation*}
\label{rmk:dual_R}
\end{remark}

\noindent \textbf{Proof.} $\Res$ is clearly linear and maps nonnegative functions to nonnegative functions. 
It maps continuous functions to continuous functions since $p^\kappa(t,x,y)$ and $\gamma$ are continuous. Thus $\Res$ is a positive linear operator mapping $C(M)\to C(M)$.
For the constant function $1: x\mapsto 1$ we have
\begin{equation}
\Res 1(x) = \int_0^\infty \dif t \int_M \Gamma(\dif y) \,p^\kappa(t,x,y)=\int_0^\infty \dif t\,\P_x(\tau_\partial > t)=\E_x[\tau_\partial].
\label{eq:Res_1}
\end{equation}
Since we are assuming that the killing rate is everywhere bounded below by $\ubar{\kappa}$, it follows that we have the uniform bound over $x \in M$,
$$ 
\E_x[\tau_\partial] \le 1/\ubar{\kappa}.
$$ 
Thus since $\Res$ is positive, it follows that $\Res$ is bounded. $\Box$ \\

\noindent We note for future reference that
\begin{equation}
\frac{1}{\bar \kappa} \le \inf_{x\in M} \Res 1(x)\le \sup_{x\in M} \Res 1(x) \le \frac{1}{\ubar\kappa}.
\label{E:R_bounded}
\end{equation}

\begin{remark}
Heuristically, the resolvent describes the average cumulative occupation measure 
of the killed diffusion over a single lifetime.
\end{remark}

Fix a probability measure $\mu\in\Pr$. We now define the fixed rebirth process with rebirth distribution $\mu$, abbreviated to FR($\mu$), to be a Markov process with c\`{a}dl\`{a}g paths, 
evolving according to the SDE \eqref{eq:SDE} between regeneration events, which occur at rate $\kappa(X^\mu_t)$. At such an event the location is drawn independently from distribution $\mu$. It can be constructed explicitly as in the construction at the end of Section \ref{subsec:assumps} in the simpler case when $X_{T_{n+1}}\sim \mu$ for each $n$. Let $(P^\mu_t)_{t\ge 0}$ denote the semigroup of this process.

Since $\kappa$ is continuous and bounded and $A$ is smooth, it can be easily shown that $(P_t^\mu)_{t \ge 0}$ is a strongly continuous semigroup on $C(M)$, so the FR($\mu$) is a Feller--Markov process, and so by the Hille--Yosida theorem (Theorem 1.7, \cite{Gross2006}) it has an infinitesimal generator $-L_\mu$ defined on a dense domain $\mathcal D(L_\mu)\subset C(M)$.

The action of the generator $-L_\mu$ on smooth functions agrees with \eqref{eq:gen_FRmu}. Since $\kappa$ is bounded and continuous, $\mathcal D(L_\mu)$ will consist of twice continuously differentiable functions on $M$, and in fact $\mathcal D(L_\mu)$ is independent of $\mu$.

Recall that $Y$ denotes the unkilled process that evolves according to the SDE $\eqref{eq:SDE}$. Since the FR($\mu$) process exhibits a natural renewal behavior, 
by conditioning on the first arrival time $\tau_\partial$, we see that
\begin{equation}
    P_t^\mu f(x) =  \int_0^t  \E_x\left [ \kappa(Y_s) \mathrm e^{-\int_0^s \kappa(Y_u)\dif u}\right ]\mu P_{t-s}^\mu f\dif s + \mathbb E_x\left [f(Y_t)  1_{\{\taud >t\}}\right ].
    \label{eq:Pmut_alternative}
\end{equation}
The second term can be expressed equivalently as in \eqref{eq:semigp_alt}.


\begin{lemma}
Given $\mu\in \Pr$, an invariant measure for the FR($\mu$) process is given by
\begin{equation} \label{E:Pimu}
\Pi(\mu)(f)=\frac{\mu\Res f}{\mu\Res 1}. 
\end{equation}
\label{lem:inv_pi}
\end{lemma}
\noindent\textbf{Proof.} 
Let $f\in \mathcal D(L_\mu)$, then it follows that $f \in \mathcal D(L^\kappa)$. We wish to show that $\mu\Res L_\mu f=0$ (Proposition 9.2 of \cite[Chapter 9]{Ethier1986}). Note by \eqref{eq:semigp_alt} that we can write
\begin{equation*}
\Res f = \int_0^\infty e^{-t\Lk}f \dif t.
\end{equation*}
Then
\begin{align*}
-\mu\Res L_\mu f &= -\mu \int_0^\infty e^{-t\Lk} L_\mu f \dif t \\
&= \mu \int_0^\infty e^{-t\Lk}\bigg (-\Lk f +\kappa \int f(y)\,\mu(\dif y) \bigg) \\
&= \mu \bigg( [e^{-t\Lk} f]_0^\infty +\mu(f)\int_0^\infty e^{-t\Lk}\kappa\dif t \bigg)
\end{align*}
where we used the backward equation (see, for instance, \cite[Chapter 1]{Ethier1986}, Proposition 1.5),
\begin{equation*}
\frac{\dif}{\dif t} (e^{-t\Lk} f)=-e^{-t\Lk}\Lk f.
\end{equation*}
Note that by Tonelli's theorem we can exchange the order of integration to find
\begin{equation*}
\int_0^\infty e^{-t\Lk} \kappa\dif t = \E\bigg[ \int_0^{\tau_\partial}\kappa(Y_s)\dif s\bigg]=\E[\xi]=1
\end{equation*}
where $\xi\sim$ Exp(1) by the definition of our killing construction.

Thus putting the terms together we have that
\begin{equation*}
-\mu\Res L_\mu f = \mu (-f +\mu(f)1)=-\mu(f)+\mu(f)=0.
\end{equation*}
Thus it follows that $\Pi(\mu)$, which is the normalized version of the measure $\mu\Res$, is an invariant probability measure for the FR($\mu$) process. $\Box$

\begin{proposition}
We have the bound 
\begin{equation*}
\|\nu_1 P_t^\mu - \nu_2 P_t^\mu\|_{1}\le 2 e^{-t\ubar{\kappa}}
\end{equation*}
for any $\nu_1,\nu_2, \mu \in \Pr$. In particular choosing $\nu_2 = \Pi(\mu)$ gives the bound
\begin{equation}
\|\nu P_t^\mu - \Pi(\mu)\|_{1}\le 2 e^{-t\ubar{\kappa}}
\label{eq:FR_contrac}
\end{equation}
for any $\nu, \mu\in \Pr$. It follows that $\Pi(\mu)$ is the unique invariant probability measure for the FR($\mu$) process.
\label{prop:exp_conv_Pmu}
\end{proposition}
\noindent\textbf{Proof.} This will follow straightforwardly from the coupling inequality, see, for instance, \cite[Section 4.1]{Roberts2004Gen}. This states that $\|\mathcal L(X)-\mathcal L(Y)\|_{1}\le 2\P (X\neq Y)$ for a coupling $\P$ of random variables $X, Y$ with laws $\mathcal L(X), \mathcal L(Y)$ respectively.

Since we are assuming the killing rate $\kappa$ is bounded below by $\ubar{\kappa}$ and the rebirth distribution is fixed, we can couple two processes started from different initial distributions at the first arrival time of a homogeneous Poisson process with rate $\ubar{\kappa}$. $\Box$

\begin{lemma}
The map $\Pi:\Pr \to \Pr$ is continuous in the topology of weak-* convergence.
\label{lem:cty_Pi}
\end{lemma}
\noindent\textbf{Proof.} 
We know that $\Res$ is a bounded linear operator on $C(M)$ by Lemma \ref{lem:Res_basic}. Thus as noted in Remark \ref{rmk:dual_R} it acts dually on the space of finite signed Borel measures, and is continuous on the dual space. So it must also be weak-* continuous. Continuity of $\mu\mapsto \Pi(\mu)$ follows. $\Box$

\begin{remark}
    In fact, since $\Pr$ is a compact metric space it follows that $\Pi:\Pr \to \Pr$ is uniformly continuous. 
    \label{rem:unif_cty_Pi}
\end{remark}

\begin{proposition}
$\mu\in\Pr$ satisfies the fixed point equation
\begin{equation}
\mu=\Pi(\mu)
\label{eq:fixed_pt}
\end{equation}
if and only if $\mu$ is quasi-stationary for the diffusion $Y$ killed at rate $\kappa$.

There exists a unique quasi-stationary distribution $\pi$ for the diffusion $Y$ killed at rate $\kappa$. Furthermore, $\pi$ has a strictly positive $C^\infty$ density with respect to the Riemannian measure, which will also be denoted by $\pi$.
\label{prop:QSD_exist}
\end{proposition}

\noindent\textbf{Proof.} Suppose $\mu$ is invariant for $L_\mu$. This means that for all smooth $f$, 
$$ \mu L_\mu f =0;$$
that is, by \eqref{E:Pimu}
$$ \frac{1}{2}\mu(\Delta f)+\mu(\nabla A\cdot \nabla f) +\mu(\kappa)\mu(f)-\mu(\kappa f)=0,$$
which is equivalent to 
$$ \mu\Lk f = \mu(\kappa) \mu(f).$$ 
Since $\mu(\kappa)>0$ this tells us that $\mu$ is a quasi-stationary distribution for $X$.

Conversely, suppose $\mu$ is quasi-stationary for $Y$. Then $$ \mu\Lk f=\lambda_0 \mu(f)$$ for all smooth $f$ and some $\lambda_0>0$. Then choosing $f\equiv 1$, we find that $\mu(\kappa)=\lambda_0$, from which it follows that $\mu L_\mu f=0$. Hence $\mu$ is stationary for $L_\mu$.

Existence and uniqueness of the quasi-stationary distribution $\pi$ follows from Theorem 1.1 of \cite{Champagnat2016a}. For each $t>0$, the transition subdensity $p^\kappa(t,x,y)$ is strictly bounded away from 0; we have that for fixed $t>0$, we can find some $\epsilon>0$ such that $p^\kappa(t,x,y) > \epsilon$, for all $x,y \in M$. Then the Assumption (A) of Theorem 1.1 of \cite{Champagnat2016a} is easily verified.


\sloppy In addition, it follows from the existence of the transition subdensity $p^\kappa(t,x,y)$ that $\pi$ must also be absolutely continuous with respect to the Riemannian measure: Since it is the quasi-stationary distribution of the diffusion $Y$, by basic properties of quasi-stationary distributions (e.g. Theorem 2.2 of \cite{Collet2013}) there exists some $\lambda_0^\kappa >0$ such that for each measurable $A \subset M$ and $t>0$,
\begin{equation*}
    \pi(A) = \P_\pi (Y_t \in A\,| \taud>t) = \int_M \pi(\dif x)\int_A  \dif y \, \gamma(y)\mathrm e^{t \lambda_0^\kappa} p^\kappa(t,x,y).
\end{equation*}
In particular this implies that $\pi$ is absolutely continuous with respect to the Riemannian measure; hence $\pi$ has a density with respect to it, which we will also denote by $\pi$. Thus the density $\pi$ satisfies for each $t>0$,
\begin{equation}
    \pi(y) = \int_M \pi(\dif x) \,\mathrm e^{t \lambda_0^\kappa} p^\kappa(t,x,y) \gamma(y),
    \label{eq:pi_p}
\end{equation}
for almost every $y \in M$. But since $p^\kappa(t,x,y)$ is positive and smooth, it follows that the density $\pi$ is continuous -- so \eqref{eq:pi_p} holds for all $y\in M$ -- and then the density $\pi$ is smooth, and also positive everywhere. $\Box$


\begin{remark}
As noted in Section \ref{sec:prelim}, in the context of Monte Carlo the quasi-stationary distribution $\pi$ will be chosen to coincide with a distribution of interest by choosing $\kappa$ according to \eqref{eq:kappa}.
\end{remark}

 \section{Deterministic Flow}
 \label{sec:flow}
\subsection{Basic properties}
We are now in a position to define the deterministic measure-valued flow that will characterise the asymptotic behavior of the normalized occupation measures $(\mu_t)_{t\ge 0}$.

Recall, as in \cite{Benaim1999}, that on a metric space $E$ a \textit{semiflow} $\Phi$ is a jointly continuous map
$$
	\Phi:\R_+ \times E\to E,
$$
$$
(t,x)\mapsto \Phi(t,x) = \Phi_t(x)
$$
such that $\Phi_0$ is the identity on $E$ and $\Phi_{t+s}=\Phi_t\circ\Phi_s$ for all $s,t \in \R_+$.

We would like to define a semiflow $\Phi$ on the space $E=\Pr$ of probability measures with the topology of weak-* convergence, which is a metric space. In particular we want $t\mapsto \Phi_t(\mu)$ to solve the measure-valued ordinary differential equation (ODE)
\begin{equation}
\dot\nu_t =-\nu_t +\Pi(\nu_t),\quad \nu_0=\mu,
\label{eq:flow}
\end{equation}
in the weak sense, meaning that for any test function $f\in C(M)$ 
$$
	\Phi_t(\mu)f=\mu f +\int_0^t \big(-\Phi_s(\mu)f+\Pi(\Phi_s(\mu))f \big)\dif s.
$$

We define such a semiflow by adapting the approach of \cite[Section 5]{Benaim2016} to our present setting.
As noted in Lemma \ref{lem:Res_basic}, the operator $\Res$ is a bounded linear operator mapping from the Banach space $C(M)$ to itself. This allows us to define, for any $t\ge 0$, the bounded linear operator $e^{t\Res}: C(M)\to C(M)$, whose dual acts on the space of finite signed Borel measures, equipped with the total-variation norm.

This allows us to define, for each $t\ge 0$, the probability measures
\begin{equation*}
\tilde \Phi_t(\mu) := \frac{\mu e^{t\Res}}{\mu \etR 1}.
\end{equation*}
The map $t\mapsto \tilde \Phi_t(\mu)$ satisfies the weak measure-valued ODE
\begin{equation*}
\dot \nu_t = \nu_t \Res - (\nu \Res 1)\nu_t=(\nu_t \Res 1)(-\nu_t +\Pi(\nu_t)),\quad \nu_0= \mu.
\label{eq:reduced_ODE}
\end{equation*}
To get a solution to our actual ODE \eqref{eq:flow} we employ a suitable time change,
imitating \cite{Blanchet2016} and \cite{Benaim2016}. Similarly to \cite{Benaim2016}, for $t\ge 0$ set
\begin{equation*}
s_\mu(t):= \int_0^t \tilde \Phi_s(\mu)\Res 1\dif s
\end{equation*}
so $\dot s_\mu (t)=\tilde \Phi_t(\mu)\Res 1= \E_{\tilde \Phi_t(\mu)}[\taud]>0$, so $t\mapsto s_\mu(t)$ is strictly increasing.
\begin{lemma}
$\dot s_\mu(t)\ge 1/\bar \kappa>0$ for all $\mu\in \Pr$ and $t>0$. Thus in particular $s_\mu(t)\to \infty$ as $t\to \infty$ for any $\mu\in\Pr$.
\end{lemma}
\noindent\textbf{Proof.} Immediate since the function $x\mapsto \E_x[\tau_\partial]$ is uniformly bounded below by $1/\bar \kappa$ since the killing rate is bounded above by $\bar \kappa$. $\Box$\\

For a fixed $\mu\in\Pr$, since $s_\mu$ is a strictly increasing, differentiable map $\R_+\to\R_+$ we can define an inverse mapping $\tau_\mu$, and compose
\begin{equation*}
\Phi_t(\mu) = \tilde \Phi_{\tau_\mu (t)}(\mu).
\end{equation*}

Recall that we have equipped the space $\Pr$ with the weak-* topology, which is a compact metric space.
\begin{proposition}
$\Phi$ is an injective semiflow on $\Pr$, and for each $\mu\in \Pr$ $t\mapsto \Phi_t(\mu)$ is the unique weak solution to \eqref{eq:flow}.
\end{proposition}
\noindent \textbf{Proof.} This is identical to the proof of \cite[Proposition 5.1]{Benaim2016}. $\Box$

\subsection{Stability of $\pi$}
Recall that we are interested in working with respect to the measure $\Gamma(\dif y)=\gamma(y)\dif y$ where $\gamma= \exp(2A)$. From Proposition \ref{prop:QSD_exist}, we see that there is a (unique) quasi-stationary distribution $\pi$ which is a fixed point of $\Pi$, which also has a density with respect to the Riemannian measure, which we also denote by $\pi$. Set
$$
\varphi:= \pi /\gamma
$$
which is a smooth function.

We now list some basic facts about $\pi$ and $\varphi$ which can all be verified through routine manipulations. There exists $\lambda_0^\kappa >0$, which describes the asymptotic rate of killing, such that $$\Lk \varphi = \lambda_0^\kappa \varphi. $$
$\varphi$ is a bounded above and bounded away from 0 since by Proposition \ref{prop:QSD_exist}, $\pi$ is both bounded above and bounded away from zero.


Writing $$\beta:= \frac{1}{\lkap}$$ it follows that 
\begin{equation}
\begin{split}
\Res \varphi &= \beta\varphi, \\
\pi \Res &= \beta \pi
\end{split}
\label{eq:Res_eigen}
\end{equation}
where this final identity holds both in terms of $\pi$ as a measure, and pointwise as density functions.

We wish now to analyse the asymptotic behavior of the semiflow $\Phi$ defined in Section \ref{sec:flow}. To do this we will derive a drift condition for a reweighted version of the Markov process.

Since $\varphi$ is bounded above and away from zero, consider the bounded linear operator $L: C(M) \to C(M)$ given by
\begin{equation*}
L= \frac{1}{\varphi}(\Res-\beta)\varphi.
\end{equation*}
Thus by exponentiation $L$ generates a Markov semigroup $(\Kt)_{t\ge 0}$. For each $t\ge 0$,
\begin{align*}
\Kt f = \frac{1}{\varphi}\exp(t(\Res-\beta))(\varphi f).
\end{align*}
We can also define the kernels $(\tilde \Kt)_{t\ge 0}$,
\begin{equation*}
\tilde \Kt := \varphi \Kt (f/\varphi) = \exp(t(\Res - \beta))f.
\end{equation*}

For the process defined by $(\Kt)_{t\ge 0}$ it can easily be seen from \eqref{eq:Res_eigen} that the measure $\pi\varphi$ given by $(\pi\varphi)(f)=\int f(x)\varphi(x)\pi(x)\dif x$ is an invariant measure. Since $\varphi$ is bounded, we see that $\pi\varphi$ is in fact a finite measure. In what proceeds, without loss of generality we rescale $\varphi$ so that $\pi \varphi(1)=1$.

We would like to show that this process is $V$-uniformly ergodic, with $V=1/\varphi$. We do this using a drift condition from Theorem 5.2 of \cite{Down1995}: for a continuous-time irreducible aperiodic Markov process with extended generator $L$, assume it satisfies for constants $b,c>0$ and petite set $C$
\begin{equation}
LV\le -cV+b1_C,
\label{eq:Vdrift_cond}
\end{equation}
then the process is $V$-uniformly ergodic. Heuristically, a petite set is a set from which the Markov process leaves with a common minorizing measure. The precise definition is carefully presented in \cite[Section 3]{Down1995}. For our present purposes, since we have a positive continuous transition subdensity $p^\kappa(t,x,y)$, it is enough to note that compact sets (and hence the entire space) are petite.

\begin{proposition}
The drift condition \eqref{eq:Vdrift_cond} holds.
\end{proposition}
\noindent\textbf{Proof.} Recall we have set $V=1/\varphi$. Then for $\beta = 1/\lambda_0^\kappa$,
\begin{equation*}
LV=\frac{1}{\varphi}(\Res-\beta)1=\frac{1}{\varphi}(\E_\cdot [\tau_\partial]-\beta)=-\beta V+V\E_\cdot [\tau_\partial]\le -\beta V + b
\end{equation*}
where $b$ is an upper bound for $V(x)\E_x[\tau_\partial]$ for all $x$, which exists since $V$ is bounded above and $\kappa$ is bounded away from 0. Note our entire space is compact, hence petite.
Thus the drift condition holds. $\Box$ 

\begin{remark}
    We note that since the entire space $M$ is compact hence petite, the drift condition \eqref{eq:Vdrift_cond} can be trivially satisfied by choosing $V\equiv 1$, the constants $b=c=1$ and $C=M$. However we have kept the choice $V=\varphi$ in the proof since this choice is suggestive of how one might generalise this to noncompact spaces.
\end{remark}

By \cite[Theorem 5.2]{Down1995} this implies $V$-uniform ergodicity: There exist constants $D$ and $0\le \rho<1$ such that for all $x\in M$:
\begin{align*}
\sup_{|g|\le V}| \Kt(x,g)-\pi\varphi(g)|\le V(x)D\rho^t. 
\end{align*}
Multiplying through by $\varphi(x)$ and relabeling $\varphi g$ as $f$, we see that
the condition $|g|\le V=1/\varphi$ is equivalent to 
$|f|\le 1$, hence proving uniform ergodicity: For any $x\in M$
\begin{equation*}
\sup_{|f|\le 1}|\tilde \Kt(x,f)-\varphi(x)\pi(f)| \le D\rho^t.
\end{equation*}
Thus we will have, for any initial distribution $\mu$,
\begin{equation}
\sup_{|f|\le 1}|\mu \tilde \Kt f - \mu(\varphi)\pi(f)|\le D\rho^t.
\label{eq:conv_tildeK}
\end{equation}

\begin{proposition}
We have convergence $\tilde \Phi_t(\mu)\to \pi$ as $t\to\infty$ in total variation distance,
uniformly in $\mu$ .
\label{prop:conv_tild_Phi}
\end{proposition}
\noindent \textbf{Proof.} 
Let $\varphi_*:= \min_{x\in M} \varphi(x)$, which is
positive, since $\varphi$ is a continuous positive
function on a compact set. We find that for any
${t\ge (\log 2D - \log \varphi_* )/\log\rho^{-1}}$
we have by \eqref{eq:conv_tildeK}
$$
	\mu\tilde \Kt 1 \ge \mu(\varphi) - D \rho^t \ge \frac{\varphi_*}{2}.
$$
Since
$$
\tilde \Phi_t(\mu) =\frac{\mu \tilde \Kt}{\mu \tilde \Kt 1},
$$
we have then for any probability measure $\mu$ and continuous $f$ with $|f|\le 1$
\begin{align*}
|\tilde\Phi_t(\mu)f-\pi(f)|&\le \left( \mu\tilde \Kt 1 \right)^{-1}\left[ \bigl|\mu\tilde \Kt f - \mu(\varphi)\pi(f) \bigr| + | \pi(f) |\cdot \bigl| \mu(\varphi)-\mu \tilde \Kt 1 \bigr|\right]\\
&\le \frac{2}{\varphi_*} \cdot \left[ 1 + |\pi(f) |\right] D \rho^t \\
&\le \frac{4D}{\varphi_*} \rho^t .
\end{align*}
$\Box$ \\

\noindent Finally, we see that this convergence carries over to the semiflow $\Phi$.
\begin{theorem}
We have convergence $\Phi_t(\mu) \to \pi$ as $t\to\infty$ uniformly in $\mu$ in total variation distance.
\end{theorem}

\noindent\textbf{Proof.} This follows from Proposition \ref{prop:conv_tild_Phi} and the fact that $\dot s_\mu(t)$ is bounded above by $1/\ubar{\kappa}$ uniformly in $\mu$ and $t$, as in Lemma \ref{lem:Res_basic}. The boundedness of this derivative ensures that its inverse $\tau_\mu$ satisfies for all $\mu \in \Pr$ and $t\ge 0$
\begin{equation*}
\tau_\mu(t) \ge \ubar{\kappa} t.
\end{equation*}$\Box$ \\

\begin{remark}
In the language of \cite{Benaim1999}, this shows that $\pi$ is a \textit{global attractor} of the semiflow $\Phi$.
\label{rem:global_attrac}
\end{remark}

\section{Asymptotic Pseudo-Trajectories}
\label{sec:APTs}
\subsection{Basic properties and definitions}
We now wish to relate the stochastic behavior of the normalised weighted empirical measures 
$$
\mu_t=\frac{r\mu_0 +\int_0^t \eta_s\delta_{X_{s-}}\dif s}{r+\int_0^t \eta_s \dif s}
$$ 
to the deterministic behavior of the flow defined in Section \ref{sec:flow}.\\

\noindent\textbf{Definition.} (\cite[Section 3]{Benaim1999}) For a metric space $(E,d)$, given a semiflow $\Phi$ on $E$, a continuous function $w:[0,\infty) \to E$ is an \textit{asymptotic pseudo-trajectory of $\Phi$} if for all $T>0$,
$$
\lim_{t\to \infty} \sup_{s\in[0,T]} d\big (w(t+s), \Phi_s(w(t))\big ) = 0.
$$

Recall that by Assumption \ref{assump:eta}, since $g$ is continuously differentiable, strictly increasing and $g(t)\to\infty$ as $t\to \infty$, it is a diffeomorphism of $\R_+ \to \R_+$. Thus we can define its inverse function $g^{-1}:\R_+ \to \R_+$, which is also continuously differentiable and satisfies $g^{-1}(t)\to\infty$ as $t\to \infty$. Let
\begin{equation*}
\zeta_t := \mu_{h(t)}
\end{equation*}
for all $t\ge 0$, where we defined $h(t):= g^{-1}(r\mathrm{e}^t-r)$. We will show that almost surely $t\mapsto \zeta_t$ is an asymptotic pseudo-trajectory of the semiflow $\Phi$ defined in Section \ref{sec:flow}.

Since 
\begin{equation}
\frac{\dif \phantom{t}}{\dif t}{\mu_t} = \alpha_t \big(-\mu_t + \delta_{X_{t-}} \big),
\label{E:DE_mu}
\end{equation}
applying the chain rule and product rule for derivatives yields
\begin{equation}
\frac{\dif \phantom{t}}{\dif t}{\zeta_t}= \big (-\zeta_t+\Pi(\zeta_t)\big )+\big (\delta_{X_{h(t-)}}-\Pi(\zeta_t)\big ).
\label{eq:DE_zeta}
\end{equation}
Looking at the first bracket, we recognize the flow from Section \ref{sec:flow}. Thus to formally show that $\zeta_t$ approximates the flow, we need to control the second bracket.
We also note here for future reference that
\begin{equation} \label{E:DE_Pit}
\begin{split}
\frac{\partial\phantom{t}}{\partial t} \Pi(\mu_t)f&= \alpha_t \bigg( -\frac{\mu_t \Res f}{\mu_t \Res 1}
+\frac{\Res f (X_{t-})}{\mu_t \Res 1}
+\frac{\mu_t \Res f}{\mu_t \Res 1}
-\frac{\mu_t \Res f \cdot \Res 1 (X_{t-})}{(\mu_t \Res 1)^2} \biggr)\\
&= \frac{\alpha_t}{\mu_t \Res 1} \bigg( \Res f(X_{t-})-\Res 1 (X_{t-})\Pi(\mu_t)f\bigg).
\end{split}
\end{equation}

We formalize this intuition by the approach of \cite[Proposition 3.5]{Benaim2002} and \cite[Lemma 5.4]{Kurtzmann2010}. 
It is proven in Theorem
3.2 of \cite{Benaim1999} that asymptotic pseudo-trajectories must be uniformly continuous. Conversely, Theorem
3.2 of \cite{Benaim1999} also tells us that a
uniformly continuous path $\zeta$ is an asymptotic
pseudo-trajectory if and only if every limit point
of the time shifts $\Theta^t\zeta$ in the
topology of uniform convergence on compact sets is itself a
trajectory of the flow.
We define $\mathcal M(M)$ to be the space of Borel signed measures on $M$, equipped with the weak-* topology, which can be metrized analogously to \eqref{eq:metric}.
Let $C(\R_+,\Pr)$ and $C(\R_+, \mathcal M(M))$ be the spaces of continuous paths mapping $\R_+$ into $\Pr$ and $\mathcal M(M)$ respectively, each equipped with the topology of uniform convergence on compact subsets of $\R_+$.
As usual, for each $t\ge 0$ we define $\Theta^t:C(\R_+,\Pr)\to C(\R_+,\Pr)$ to be the shift map given by 
$$ 
\bigl[\Theta^t \zeta\bigr]_s=\zeta_{t+s}, \quad s \ge 0.
$$
Defining the \textit{retraction} $\hat \Phi: C(\R_+, \mathcal P(M)) \to C(\R_+, \mathcal P(M))$ as in \cite{Benaim1999} and \cite[Proposition 3.5]{Benaim2002} by
\begin{equation*}
    \hat\Phi (\zeta)(s) = \Phi_s(\zeta_0), \quad s\ge 0,
\end{equation*}
showing that $\zeta$ is an asymptotic pseudo-trajectory of $\Phi$ is then equivalent
to showing that the limit points of
$\{\Theta^t \zeta\}_{t\ge 0}$ are fixed points of the retraction $\hat \Phi$.

We also define in analogue to \cite{Benaim2002} the operator
$L_F: C(\R_+, \Pr)\to C(\R_+, \mathcal M(M))$ by
$$
	L_F(\nu)(s)=\nu(0)+\int_0^s F(\nu(u))\dif u, \quad s \ge 0,
$$ 
with $F(\nu):=-\nu+\Pi(\nu)$, which will be used in the subsequent proof.

Define the collection of signed measures 
$$
\epsilon_t(s):=\int_t^{t+s} \big (\delta_{X_{h(u-)}} - \Pi(\zeta_u)\big)\dif u
$$ 
for all $t,s\ge 0$. We note that for each $t\ge 0$, $\epsilon_t(\cdot)\in C(\R_+, \mathcal M(M))$.

\begin{theorem}
Suppose $\zeta: \R_+ \to \Pr$ is a continuous path as described above.
Then $\zeta$ is an asymptotic pseudo-trajectory for $\Phi$ if and only if we have the following condition:

For any $T>0$ and smooth $f\in C(M)$ 
\begin{equation}
\lim_{t\to \infty} \sup_{s\in[0,T]}|\epsilon_t(s)f|=0.
\label{eq:eps_conv}
\end{equation}
\label{thm:APT_charac}
\end{theorem}
\begin{remark}
Theorem 3.2 of \cite{Benaim1999} assumes relative compactness of the image of the path $\zeta$. By Prokhorov's theorem (see \cite[Theorem 11.5.4]{Dudley2002}), relative compactness in $\Pr$ is equivalent to tightness, which trivially holds in our present compact setting.
\end{remark}

\noindent\textbf{Proof.} Given continuous $f$ with $\|f\|_\infty\le 1$, we have 
$$ 
|\zeta_{t+s} f -\zeta_t f|\le 2 |s|
$$ 
since 
$$
\frac{\dif (\zeta_t f)}{\dif t}=-\zeta_t f + f(X_{h(t-)}). 
$$ 
Hence we have uniform continuity.

Suppose the condition \eqref{eq:eps_conv} holds for any $T>0$ and smooth $f\in C(M)$. This says that $\epsilon_t(\cdot)$ converges to 0 in $C(\R_+, \mathcal M(M))$.
Consider the family $\{\Theta^t \zeta \}_{t\ge 0}$. Suppose $\tilde \zeta$ is a limit point of this family in $C(\R_+,\Pr)$.
Analogously to \cite{Benaim2002}, we can use \eqref{eq:DE_zeta} and the definition of $\epsilon_t$ to write 
\begin{equation}
\Theta^t \zeta = L_F(\Theta^t \zeta ) + \epsilon_t(\cdot).
\label{eq:Theta_LF}
\end{equation}
Since we are assuming precisely \eqref{eq:eps_conv}, and since $L_F$ is continuous, taking $t\to\infty$ in \eqref{eq:Theta_LF} along the appropriate subsequence we obtain $\tilde{\zeta} = L_F(\tilde \zeta)$. This shows that the limit path $\tilde{\zeta}$ is a fixed point of $L_F$, that is,
\begin{equation*}
    \tilde \zeta_s = \tilde \zeta_0 + \int_0^s F(\tilde \zeta_u)\dif u, \quad s \ge 0.
\end{equation*}
By uniqueness of the flow this implies that $\tilde \zeta_s = \Phi_s (\tilde \zeta_0)$ for all $s \ge 0$, that is, $\tilde \zeta$ is a fixed point of the retraction $\hat \Phi$. This concludes the proof of the sufficiency of condition \eqref{eq:eps_conv}.

Conversely, suppose $\zeta$ is an asymptotic pseudo-trajectory for $\Phi$. 
By definition, this means that for each $T>0$, as $t\to \infty$,
\begin{equation}
    \sup_{s \in [0,T]} d_w\left (\zeta_{t+s}, \Phi_s(\zeta_t) \right) \to 0.
    \label{eq:zeta_APT}
\end{equation}
By the representation \eqref{eq:Theta_LF} we would like to show that for each $T>0$,
\begin{equation*}
    \sup_{s\in [0,T]} d_w\left(\Theta^t \zeta (s), L_F(\Theta^t \zeta)(s) \right) =\sup_{s\in [0,T]} d_w\left(\zeta_{t+s},\, \zeta_t+ \int_0^s F(\zeta_{t+u})\dif u\right)\to 0
\end{equation*}
By \eqref{eq:zeta_APT}, it is sufficient to control $d_w(L_F(\Theta^t \zeta_s), \Phi_s(\zeta_t))$ uniformly over $s \in [0,T]$. 

By the definition of the flow $\Phi$, this will follow if we can control
\begin{equation}
    \sup_{s\in [0,T]} d_w\left(\int_0^s F(\zeta_{t+u})\dif u , \int_0^s F(\Phi_u(\zeta_t))\dif u\right).
    \label{eq:LF_Phi_disc}
\end{equation}
Since $F$ is uniformly continuous (since $\Pi: \Pr\to\Pr$ is uniformly continuous; Remark \ref{rem:unif_cty_Pi}), \eqref{eq:LF_Phi_disc} can be made arbitrarily small by choosing $t$ sufficiently large because of \eqref{eq:zeta_APT}.

Thus by \eqref{eq:Theta_LF} we conclude that $\epsilon_t \to 0$ in $C(\R_+, \mathcal M(M))$, that is, \eqref{eq:eps_conv} must hold for any $T>0$ and smooth $f\in C(M)$ as required. $\Box$\\

\noindent The goal now is to establish the requirements of Theorem \ref{thm:APT_charac} almost surely. So we need to establish \eqref{eq:eps_conv}.

\subsection{Poisson Equation}
To show that \eqref{eq:eps_conv} holds, inspired by \cite[Section 5.2]{Benaim2002} and \cite[Section 6.2]{Benaim2016}, we will measure the discrepancy between the cumulative
occupancy and the quasi-stationary distribution by a solution to the Poisson equation
\begin{equation}
-f+\Pi(\mu)f=-L_\mu g.
\label{eq:Poisson_eq}
\end{equation}

\noindent Fix any $\mu\in\Pr$. 
Define for $f\in C(M)$, 
$$Q_\mu f := \int_0^\infty \left (P_t^\mu f -\Pi(\mu)f \right)\dif t.$$
\begin{lemma}
For any $\mu\in\Pr$, $Q_\mu$ is a bounded linear operator mapping $C(M)$ to itself, and $$\|Q_\mu f\|_\infty \le  \frac{2\|f\|_\infty}{\ubar{\kappa}}$$ for any $\mu\in \Pr$ and $f\in C(M)$.
\label{lem:Q_mu}
\end{lemma}

\noindent\textbf{Proof.} $Q_\mu f$ is continuous since $f$ is continuous and $(P_t^\mu)$ is Feller. From \eqref{eq:FR_contrac}, we have that 
\begin{equation*}
|Q_\mu f(x)| \le \int_0^\infty \|P_t^\mu(x)  -\Pi(\mu)\|_{1} \|f\|_\infty\dif t \le 2\|f\|_\infty \int_0^\infty e^{-t\ubar{\kappa}}\dif t=\frac {2\|f\|_\infty}{\ubar{\kappa}}
\end{equation*}
uniformly over $x\in M$ and $\mu\in \Pr$. $\Box$ \\

\noindent Recall that $-L_\mu$ is the generator of the FR($\mu$) process. By the basic properties of infinitesimal generators (see, for instance, \cite{Gross2006}), we have that for any $t>0$, $f\in C(M)$, that $\int_0^t P_s^\mu f \dif s \in \mathcal D(L_\mu)$, and $(-L_\mu) \int_0^t P_s^\mu f \dif s= P_t^\mu f - f$. Then 
\begin{equation}
P_t^\mu f(x)=f(x)+(-L_\mu) \int_0^t P_s^\mu f(x)\dif s .
\label{eq:Pmu_Amu}
\end{equation}

\begin{proposition}
    Given $f\in C(M)$, $Q_\mu f \in \mathcal D(L_\mu)$, and $g=Q_\mu f$ solves the Poisson equation \eqref{eq:Poisson_eq}.
    \label{prop:Q_mu}
\end{proposition}

\noindent\textbf{Proof.} Without loss of generality, assume $\Pi(\mu) f =0$; just replace $f$ with $f-\Pi(\mu)f$. Then taking $t\rightarrow \infty$ in \eqref{eq:Pmu_Amu} and using Proposition \ref{prop:exp_conv_Pmu} we see by closure of $-L_\mu$ that $Q_\mu \in \mathcal D(L_\mu)$, and that $\Pi(\mu) f=f-L_\mu Q_\mu f$. Thus \eqref{eq:Poisson_eq} is satisfied by $g=Q_\mu f$. $\Box$

\begin{remark}
By basic properties of generators of Feller processes (\textit{cf.}\ \cite[Section VII.1]{RevuzYor91}), we similarly have that for $f\in \mathcal D(L_\mu)$,
\begin{equation}
-f + \Pi(\mu) f = -Q_\mu L_\mu f.
\label{eq:Poisson_eq_rev}
\end{equation}
\end{remark}

\begin{remark}
Note that we can think of the operator 
$$
	K_\mu:f\mapsto  f-\Pi(\mu)f
$$ 
as a projection operator since it is linear and idempotent ($K_\mu^2=K_\mu$).
\end{remark}

\subsection{Bounding the discrepancy}
Using our solution to the Poisson equation, we now rewrite the discrepancy term, and decompose it using It\^{o}'s formula. Making the change of variables $u\leftarrow h(u)$ we write
\begin{align*}
\epsilon_t(s)f&=\int_t^{t+s} (f(X_{h(u-)}) -\Pi(\zeta_u) f )\dif u \\
&= \int_{h(t)}^{h(t+s)}(-L_{\mu_u})Q_{\mu_u}f(X_{u-})\,\alpha_u{\dif u}.
\end{align*}

Consider now a time-dependent function $f: \R\times M \to \R$, written as $f_s(x)$. For simplicity we
will take the notation $\nabla$ and $\nabla^2$ to
refer to the gradient with respect to the coordinates of $M$, and $f'_s$ to be the time derivative $\partial f_t/\partial t \bigm|_{t=s}$.
We now apply It\^o's formula for general semimartingales (as formulated, for instance, as Theorem 14.2.4 of \cite{Cohen2015}), taking advantage of the fact
that all quadratic covariation terms are 0 for the
Brownian motion that is driving our process $X_t$:
\begin{multline*}
f_t(X_t)-f_0(X_0)=\int_0^t \nabla f_u(X_{u-})\cdot \dif X_u +\frac{1}{2}\int_0^t \nabla^2 f_u(X_{u-})\dif u+\int_0^t f'_u (X_{u-})\dif u\\ 
+\sum_{0<u\le t}\big ( f_u(X_u)-f_u(X_{u-}) -\nabla f_u(X_{u-})\cdot \Delta X_u \big ).
\end{multline*}
Using the formula \eqref{eq:gen_FRmu} for $L_\mu$ we find
\begin{multline*}
f_t(X_t)-f_0(X_0)=\int_0^t (-L_{\mu_u}) f_u(X_{u-})\dif u + \int_0^t \nabla f_u(X_{u-})\cdot \dif W_u+\int_0^t f'_u(X_{u-})\dif u\\
+\sum_{0<u\le t}\big ( f_u(X_u)-f_u(X_{u-}) \big)-\int_0^t \kappa(X_{u-})\int \big (f_u(y)-f_u(X_{u-}) \big ) \mu_u(\dif y)\dif u,
\end{multline*}
where $W_u$ is a Wiener process on $M$.

We apply this formula now to the function
$f_s (x) = Q_{\mu_s} f (x) \alpha_s$. 
By Proposition \ref{prop:Q_mu}, $Q_{\mu_s}f \in \mathcal D(L_{\mu_s})$, and so in particular $Q_{\mu_s}f$ is twice continuously differentiable.
Thus $f_s(x)$ is indeed twice continuously differentiable with respect to $x$.  We rearrange the terms to obtain
\begin{equation}
\begin{split}
\epsilon_s(t)f &= \int_{h(t)}^{h(t+s)}(-L_{\mu_u}) Q_{\mu_u}f(X_{u-})\,\alpha_u{\dif u}\\
&=\epsilon_s^{(1)}(t)f+\epsilon_s^{(2)}(t)f+\epsilon_s^{(3)}(t)f+\epsilon_s^{(4)}(t)f+\epsilon_s^{(5)}(t)f,
\label{eq:eps_decomp}
\end{split}
\end{equation}
where
\begin{align*}
\epsilon_t^{(1)}(s)f&= {Q_{\mu_{h(t+s)}}f(X_{h((t+s)-)})} \, \alpha_{h(t+s)} -{Q_{\mu_{h(t)}}f(X_{h(t-)})}\,\alpha_{h(t)},\\
\epsilon_t^{(2)}(s)f&=\int_{h(t)}^{h(t+s)}{Q_{\mu_u}f(X_{u-})} \frac{\dif \alpha_u}{\dif u} \dif u,\\
\epsilon_t^{(3)}(s)f&=-\int_{h(t)}^{h(t+s)} {\frac{\partial}{\partial u}\bigg( Q_{\mu_u} f\bigg)(X_{u-})}\,\alpha_u\dif u, \\
\epsilon_t^{(4)}(s)f&= -N^f_{h(t+s)}+N^f_{h(t)},\\
\epsilon_t^{(5)}(s)f&=-J^f_{h(t+s)}+J^f_{h(t)},
\end{align*}
and $N$ and $J$ are the local martingales
\begin{align*}
N^f_t &:= \int_{0}^{t}{\nabla Q_{\mu_u}f(X_{u-})}\,\alpha_u\cdot \dif W_u ,\\
J^f_t &:= \sum_{0<u\le t}\bigg({Q_{\mu_u}f(X_u)}\alpha_u -{Q_{\mu_u}f(X_{u-})}\alpha_u \bigg) \\
&\quad-\int_{0}^{t} \dif u\,\kappa(X_{u-})\int \mu_u(\dif y)\,\bigg({Q_{\mu_u}f(y)}\alpha_u-{Q_{\mu_u}f(X_{u-})}\alpha_u \bigg).
\end{align*}

\begin{theorem}
The conditions of Theorem \ref{thm:APT_charac} hold almost surely. This implies that $t\mapsto \zeta_t$ is almost surely an asymptotic pseudo-trajectory for $\Phi$.
\end{theorem}

\noindent\textbf{Proof.} 
To establish \eqref{eq:eps_conv} we will use the decomposition \eqref{eq:eps_decomp} and consider the five error terms individually.

\subsubsection{$\epsilon_t^{(1)}(s)f$ and $\epsilon_t^{(2)}(s)f$}
Using the bound from Lemma \ref{lem:Q_mu}, we see that
\begin{align*}
|\epsilon_t^{(1)}(s)f| &\le \|Q_{\mu_{h(t+s)}}f\|_\infty \alpha_{h(t+s)} + \|Q_{\mu_{h(t)}}f\|_\infty \alpha_{h(t)} \\
&\le 4 \ubar{\kappa}^{-1}\|f\|_\infty (\alpha_{h(t+s)}- \alpha_{h(t)}),\\
|\epsilon_t^{(2)}(s)f|&\le \int_{h(t)}^{h(t+s)}{\|Q_{\mu_u}f\|_\infty} \frac{\dif \alpha_u}{\dif u}\dif u \\
 & = 2\ubar{\kappa}^{-1}\|f\|_\infty (\alpha_{h(t+s)}- \alpha_{h(t)}).
\end{align*}
Since $h(t)\to \infty$ as $t\to \infty$ and $\alpha_t\to 0$ as $t\to\infty$ by Assumption \ref{assump:eta}, we see that both terms decay to 0 as $t\to\infty$.

\subsubsection{$\epsilon_t^{(3)}(s)f$}
By the same argument as in the proof of Lemma 5.5 of \cite{Benaim2002} we have 
\begin{equation} \label{E:Qmu}
\frac{\partial}{\partial t} Q_{\mu_t} =-\bigg ( \frac{\partial}{\partial t} K_{\mu_t}+Q_{\mu_t}\frac{\partial}{\partial t}(-L_{\mu_t})\bigg )Q_{\mu_t},
\end{equation}
where $K_\mu f=f-\Pi(\mu) f$.
While our $-L_\mu$ is not the same as their operator $A_\mu$, the same proof holds, since it relies only upon the Poisson
equation \eqref{eq:Poisson_eq} and \eqref{eq:Poisson_eq_rev}.

Applying \eqref{E:DE_mu} to the definition \eqref{eq:gen_FRmu} of $L_{\mu_t}$ we obtain 
$$
-\frac{\partial}{\partial t} L_{\mu_t} f (x) = {\kappa(x)}{\alpha_t}\Bigl(f(X_{t-})-\mu_t(f)\Bigr).
$$
From
\eqref{E:DE_Pit}
we obtain an upper bound
\begin{align*}
  \left\| \alpha_t^{-1} \frac{\partial K_{\mu_t}}{\partial t} Q_{\mu_t} f   \right\|_\infty &\le
\left |\frac{1}{\mu_t \Res 1}\right | \left(\left\| \Res Q_{\mu_t} f (X_{t-}) \right\|_\infty +\left | \Pi(\mu_t) \Res 1 (X_{t-})\right |
\right),
\end{align*}
so by Lemma \ref{lem:Q_mu}
and the bounds \eqref{E:R_bounded} we see that
there is a constant $C$ (depending on the 
upper and lower bounds on $\kappa$) such that
$$
	\sup_{x\in M} \Biggl|\frac{\partial Q_{\mu_t}}{\partial t} f(x)\Biggr|
    \le C \|f\|_\infty \alpha_t .
$$
This implies that
$$
|\epsilon_t^{(3)}(s)f| \le C \|f\|_\infty \int_{h(t)}^{\infty}\alpha_u^2 \dif u.
$$
Since the definite integral from 0 to $\infty$ is finite by Assumption \ref{assump:eta}, it follows immediately that $\epsilon_t^{(3)}(s)f$ decays to 0 as $t\to \infty$.

\subsubsection{$\epsilon_t^{(4)}(s)f$}
\label{sub_sec:martingale_error}
We now turn to the first of the two martingale terms. Similarly to  Proposition 5.3 of \cite{Kurtzmann2010}, our goal is to control the quadratic variation and then apply the Burkholder--Davis--Gundy inequality.

The quadratic variation of the martingale $N_t^f$ is bounded by
\begin{equation*}
\int_0^t {\|\nabla Q_{\mu_s} f(X_s)\|_\infty^2}\,{\alpha_u^2}\dif s\, .
\end{equation*}
We can bound this by means of the inequality
\begin{equation}
\|\nabla P_t^\mu f\|_\infty \le \frac{C_1 \|f\|_\infty}{\sqrt{t}},\quad \text{for }t\text{ sufficiently small,}
\label{eq:nabla_P_bd}
\end{equation}
for some constant $C_1$. Such an inequality was shown to hold in the proof of Lemma 5.1 of \cite{Benaim2002} for diffusive processes on compact manifolds without jumps. To see that this inequality also holds in the present setting, consider the expression for the semigroup \eqref{eq:Pmut_alternative}:
\begin{equation*}
    P_t^\mu f(x) =  \int_0^t  \E_x\left [ \kappa(Y_s)  1_{\{\taud >s\}}\right ]\mu P_{t-s}^\mu f\dif s + \mathbb E_x\left [f(Y_t)  1_{\{\taud >t\}}\right ].
\end{equation*}
With this representation, it suffices to show that an inequality of the form \eqref{eq:nabla_P_bd} for the killed semigroup given in Lemma \ref{lemma:subtrans}. This will ensure that our resulting bound \eqref{eq:nabla_P_bd} is uniform over $\mu$. 

By \eqref{eq:p^kappa}, since $|g(t,x,y)|\le 1$ uniformly, and since we know that such an inequality holds for diffusive processes on compact manifolds without jumps, it suffices to show that
\begin{equation}
    \|\nabla_x g(t,\cdot,y)\|_\infty \le \frac{C}{\sqrt t}
    \label{eq:nabla_P_desired}
\end{equation}
for some constant $C$ uniformly over $y$, for all $t$ sufficiently small. By the Girsanov--Cameron--Martin formula \cite{Elworthy1982} we can write 
\begin{equation*}
    g(t,x,y) = \mathbb Q_t \left[\mathrm e^{-\int_0^t \kappa(Y_s)\dif s}\, G(Y) \right]
\end{equation*}
where for $s \in [0,t]$, $Y_s = (1-s/t) x + (s/t) y + \hat W_s$, where under $\mathbb Q_t$, $\hat W$ is a standard Brownian bridge with $\hat W_0 = \hat W_t = 0$. $G(Y)$ is given by
\begin{equation*}
    G(Y) = \exp\left (A(y)-A(x)-\frac{1}{2}\int_0^t \Delta A(Y_s)\dif s - \frac{1}{2}\int_0^t \|\nabla A(Y_s)\|^2 \dif s \right ),
\end{equation*}
since our drift is assumed to be of a gradient form $\nabla A$. Having now expressed $Y$ with explicit dependence upon the initial position $x$, we can calculate $\nabla_x g(t,\cdot, y)$. 
Since the functions $\kappa$ and $A$ are assumed smooth (hence they and their derivatives are bounded uniformly), we can conclude that a bound of the form \eqref{eq:nabla_P_desired} holds. In fact, the bound we obtain is actually uniform over small $t$, and does not blow up as $t\to 0$.

Armed with \eqref{eq:nabla_P_bd}, we then have 
\begin{equation}
\int_0^\infty \|\nabla P_t^\mu f\|_\infty \dif t \le C_2 \|f\|_\infty,
\label{eq:C_3_bound_nablaP}
\end{equation}
for some constant $C_2$, obtained by considering separately the integrals over $(0,t_0)$ and $(t_0,\infty)$.
Applying \eqref{eq:nabla_P_bd} bounds the former, while the semigroup property allows the latter piece to be bounded by
\begin{align*}
\int_{t_0}^\infty \|\nabla P_t^\mu f\|_\infty \dif t &\le \int_0^\infty \|\nabla P_{t_0}^\mu (P_t^\mu (K_\mu  f)) \|_\infty \dif t \\
&\le \frac{C_1}{\sqrt{t_0}}\int_0^\infty \| P_t^\mu (K_\mu f)\|_\infty \dif t.
\end{align*}
Here we replaced $f$ with $K_\mu f= f -\Pi(\mu) f$ since $\Pi(\mu)f$ is a constant so $\nabla \Pi(\mu)f=0$. This final term can be bounded just as in Lemma \ref{lem:Q_mu}.

Since $\nabla Q_\mu f = \int_0^\infty \nabla P_t^\mu \dif t$, \eqref{eq:C_3_bound_nablaP} immediately implies a universal bound
\begin{equation*}
\|\nabla Q_{\mu_u}f\|_\infty \le C_3 \|f\|_\infty.
\end{equation*}
for smooth $f$.
This gives us a bound for the quadratic variation of 
\begin{equation*}
 {C_3^2 \|f\|_\infty^2} \int_0^t{\alpha_s^2}\dif s.
\end{equation*}
The Burkholder--Davis--Gundy inequality then implies the existence of a constant $C_4$ such that for any $\delta>0$
\begin{equation*}
\P_x \bigg ( \sup_{s\in [0,T]} |\epsilon_t^{(4)}(s)f|\ge \delta \bigg)\le \frac{C_4\|f\|^2_\infty}{\delta^2} \int_{h(t)}^\infty \alpha_s^2 \dif s.
\end{equation*}

Using \eqref{eq:gamma_int_sum} from Assumption \ref{assump:eta},  it follows from the Borel--Cantelli lemma that almost surely, for any $\delta>0$ and $T>0$,
\begin{equation}
\limsup_{n\to\infty} \sup_{s\in[0,T]}|\epsilon_{n}^{(4)}(s)f| \le \delta.
\label{eq:eps4_conv_n}
\end{equation}
Thus for any $\delta>0$ we may find $N$ sufficiently large so that for any $t\ge N -T$,
$$
	\sup_{s\in [0,T]}\bigl|\epsilon_t^{(4)}(s)f \bigr| \le \sup_{s\in [0,1]}\bigl|\epsilon_{\lfloor t \rfloor}^{(4)}(s)f \bigr| +
    \sup_{s\in [0,T+1]}\bigl|\epsilon_{\lfloor t \rfloor}^{(4)}(s)f \bigr| +
    \sup_{s\in [0,1]}\bigl|\epsilon_{\lfloor t +T\rfloor}^{(4)}(s)f \bigr| \le \delta .
$$
Since $\delta>0$ is arbitrary, it follows that almost surely for every $T$
\begin{equation*}
\lim_{t\to\infty}\sup_{s\in[0,T]}|\epsilon_t^{(4)}(s)f| =0.
\end{equation*}

\subsubsection{$\epsilon_t^{(5)}(s)f$}
The final error term is a jump martingale term. We use the same approach as above for $\epsilon_t^{(4)}(s)f$. For this jump martingale the quadratic variation is 
\begin{equation*}
\sum_{0<u\le t} \bigg( {Q_{\mu_u}f(X_u)\alpha_u-Q_{\mu_u}f(X_{u-})\alpha_u}\bigg)^2,
\end{equation*}
where the sum is over jump-points $u\in (0,t]$.
The squared jump at time $u$ is bounded by
\begin{equation*}
16\|f\|_\infty^2 \ubar{\kappa}^{-2} \alpha_u^2 ,
\end{equation*}
making use of Lemma \ref{lem:Q_mu}. The expected
quadratic variation is then the expectation of
the predictable variation, which is bounded by
\begin{equation}
\E\bigg[ \int_{0}^{t}{\kappa(X_u)}{\alpha_u^2}\dif u\bigg]\le \bar{\kappa} \int_0^t \alpha_u^2 \dif u.
\label{eq:E_kappa}
\end{equation}
Thus the total quadratic variation is bounded, and as before we conclude that almost surely
\begin{equation*}
\lim_{t\to\infty}\sup_{s\in[0,T]}|\epsilon_t^{(5)}(s)f| =0.
\end{equation*}

\subsubsection{Concluding the proof} We have shown
that the five discrepancy terms all converge to 0 uniformly on compact sets as $t\to\infty$ almost surely. Thus condition \eqref{eq:eps_conv} holds almost surely and the proof is complete. $\Box$ \\

\subsection{Proof of Corollary \ref{cor:conv}} This follows from Theorem \ref{thm:intro_apt} and Remark \ref{rem:global_attrac}, since limit sets of asymptotic pseudo-trajectories are \textit{attractor free sets}, and so will be contained in the global attractor of $\Phi$, which is $\{\pi\}$. These relationships are spelled out in \cite[Section 5]{Benaim1999}. $\Box$

\section{Some comments on noncompact manifolds}
\label{sec:extens}
In this work we have restricted our analysis to compact state spaces. Practically speaking, though, QSMC methods such as ScaLE \cite{Pollock2016} and ReScaLE as described here are applicable on noncompact state spaces such as $\Rd$.

The key difficulties in extending this present work to the setting of 
noncompact state spaces are:
\begin{enumerate}
\item establishing almost sure tightness of the occupation measures, and
\item arguing that $\E[\int_0^t \kappa(X_s)\dif s]=O(t)$ as $t\to \infty$.
\end{enumerate}
We have also implicitly assumed there is a lower bound on the Ricci curvature, for the local bounds on the growth of the gradient of the Brownian transition kernel used in section \ref{sub_sec:martingale_error}.

To establish tightness it might be helpful to consider the discrete skeleton at the regeneration times, $(\mu_{T_n})_{n\in \mathbb N}$, as a \textit{measure-valued P\'{o}lya urn process}, as introduced in the recent work \cite{Mailler2017}.

The second point --- bounding $\E[\int_0^t \kappa(X_s)\dif s]$ --- is necessary in controlling the variance of the jump martingale in \eqref{eq:E_kappa}.

\cite{Kurtzmann2010} works in a noncompact setting and considers a related problem, namely studying self-interacting diffusions where the occupation measures $(\mu_t)$ come into play through the drift of the diffusion. However the techniques utilised --- such as the construction of explicit Lyapunov functions --- are not immediately applicable. Our choice of the metric \eqref{eq:metric} also needs to be modified, since --- contrary to what is stated in \cite[Section 2.1.1]{Kurtzmann2010}  ---
the space of bounded continuous functions is in general not separable on a noncompact state space.

\section*{Acknowledgements}
Research of A. Q. Wang is supported by EPSRC OxWaSP CDT through grant EP/L016710/1. Research of G. O. Roberts is supported by EPSRC grants EP/K034154/1, EP/K014463/1, EP/D002060/1. We would particularly like to thank the anonymous referees whose comments have substantially improved the paper.

\bibliography{compact_rescale}
\bibliographystyle{plain}

\end{document}